\documentclass[final]{elsart5p}
\journal{Comp. Meth. Applied Mech. Engineering, accepted and published as doi:10.1016/j.cma.2006.10.019}

\usepackage{ifpdf}
\usepackage{amssymb,latexsym,lineno}
\ifpdf
\usepackage[%
  pdftitle={Observations on degenerate saddle point problems},%
  pdfauthor={Andrew Knyazev},%
  pdfsubject={Degenerate saddle point problems},%
  pdfkeywords={Wellposedness,  
mixed, 
symmetric, 
saddle point, 
Lagrange multiplier, 
Ladygenskaya--Babuska--Brezzi (LBB) condition, 
inf--sup condition, 
coercivity, 
minimum gap between subspaces},%
  pdfstartview=FitH,%
  bookmarks=true,%
  bookmarksopen=true,%
  breaklinks=true,%
  colorlinks=true,%
  linkcolor=blue,anchorcolor=blue,%
  citecolor=blue,filecolor=blue,%
  menucolor=blue,pagecolor=blue,%
  urlcolor=blue]{hyperref}
\else
\usepackage[%
  breaklinks=true,%
  colorlinks=true,%
  linkcolor=blue,anchorcolor=blue,%
  citecolor=blue,filecolor=blue,%
  menucolor=blue,pagecolor=blue,%
  urlcolor=blue]{hyperref}
\fi

\usepackage[square,comma,numbers,sort&compress]{natbib} 

\newtheorem{mylemma}{Lemma}[section]
\newtheorem{mytheorem}{Theorem}[section]

\def\proof{{{\sl Proof}. }}
\def\endproof{$\Box$}

\def\R{\rm R}

\def\e{\epsilon}
\def\s{\sigma}
\def\k{\kappa}

\def\ia{{\bf R} (A)}

\def\P{{\bf P}}
\def\Pp{\P^\perp}
\def\D{{\bf D}}
\def\Dp{{\bf D}^\perp}
\def\Im{{\bf R}}
\def\Ker{{\bf N}}
\def\N{{\bf N}}
\def\ka{{\bf N} (A)}
\def\H{{\bf H}}
\def\V{{\bf V}}
\def\F{{\bf F}}

\def\R{\rm R}

\def\cal{{}}

\let\NoHyper\relax 

\begin{document}
\begin{frontmatter}

\title{Observations on degenerate saddle point problems}
\author{Andrew V. Knyazev
}
\address{Department of Mathematical Sciences\\ 
University of Colorado at Denver
and Health Sciences Center\\
P.O. Box 173364, Campus Box 170, Denver, CO 80217-3364}

\ead{andrew.knyazev[AT]cudenver.edu}
\ead[url]{http://math.cudenver.edu/ $\tilde{}$ aknyazev/}
\thanks
{Partially supported by the
National Science Foundation award DMS-0612751.}

\begin{abstract}
We investigate degenerate saddle point problems, 
which can be viewed as limit cases of 
standard mixed formulations of symmetric problems with large
jumps in coefficients. We prove that they are well-posed 
in a standard norm despite the degeneracy. 
By wellposedness we mean a stable dependence of the solution 
on the right-hand side. 
A known approach of splitting the 
saddle point problem into separate equations for the 
primary unknown and for the Lagrange multiplier is used. 
We revisit the traditional 
Ladygenskaya--Babu\v{s}ka--Brezzi (LBB) or inf--sup  condition 
as well as the standard coercivity condition, 
and analyze how they are affected by the degeneracy of the 
corresponding bilinear forms. We suggest and discuss 
generalized conditions that cover the degenerate case. 
The LBB or inf--sup  condition is 
necessary and sufficient for wellposedness of the 
problem with respect to the Lagrange multiplier 
under some assumptions. 
The generalized coercivity condition 
is necessary and sufficient for wellposedness of the 
problem with respect to the primary unknown  
under some other assumptions. 
We connect the generalized coercivity condition 
to the positiveness of the minimum gap of relevant subspaces, 
and propose several equivalent expressions 
for the minimum gap.  
Our results provide a foundation for research 
on uniform wellposedness 
of mixed formulations of symmetric problems with large
jumps in coefficients in a
standard norm, independent of the jumps.  
Such problems appear, e.g.,\ 
in numerical simulations of composite materials made of 
components with contrasting properties. 
\end{abstract}
\end{frontmatter}
\section{Introduction}
Degenerate saddle point problems, e.g.,\   
can be viewed as limit cases 
of mixed formulations of symmetric problems with large
jumps in coefficients, corresponding to an infinite jump.   
We prove that the degeneracy does not affect the wellposedness  
in a standard norm under some natural assumptions, using 
ideas that are initiated by \cite{bk90,bkk91,k92,bk94,bk95,bkp02,kw03}. 
By wellposedness, contrary to illposedness, 
we mean a stable dependence of the solution 
on the right-hand side. 
Results of this paper provide a foundation for research 
on uniform wellposedness 
of mixed formulations of symmetric problems with large
jumps in coefficients in a
standard norm, independent of the jumps.  

The necessary and sufficient condition, e.g.,\ \cite{b74,bf},
of the standard wellposedness of an operator equation 
with an arbitrary right--hand side   
is the existence of a bounded inverse of the operator. 
We argue that in some practical cases 
the equation is degenerate, i.e. 
the inverse operator does not exist.  
Assuming that the right--hand side is in the operator range, 
a solution exists, but is not unique. To make the 
solution unique we factor out the operator null--space.  
This leads to a natural generalization, 
where boundedness of the pseudoinverse of the operator is used 
as the necessary and sufficient condition of
wellposedness of a degenerate operator equation, 
by analogy with \cite{kw03,k03}. 
 
With this idea in mind, we revisit 
necessary and sufficient conditions of wellposedness 
of an abstract mixed problem. 
In the symmetric case we consider here, 
the mixed problem can be interpreted as 
a variational saddle point problem. 
For generalized  saddle point problems we refer the reader, 
e.g.,\ to \cite{MR2002909}. 

We start in Section \ref{s2} 
with a standard abstract symmetric mixed problem as in \cite{b74,bf}.
By analogy with \cite{k92,MR1971217}, 
we split the saddle point problem into two equations, 
for the primary unknown and for the Lagrange multiplier. 
This split is somewhat implicit in \cite{b74,bf}. 
The equation for the primary unknown is self-consistent, since here
we eliminate the Lagrange multiplier from the mixed system using an orthogonal projector. 

Following, e.g.,\ \cite{bf},
we discuss the traditional necessary and sufficient conditions
of wellposedness, namely,  
the Ladygenskaya--Babu\v{s}ka--Brezzi (LBB) or 
\emph{inf-sup} condition and the coercivity condition. 
The LBB or inf-sup condition, considered in Section \ref{sis}, 
is necessary and sufficient for a
stable dependence of the Lagrange multiplier 
on an arbitrary right-hand side. 

We review the traditional point of view that the coercivity condition is a
necessary and sufficient condition of wellposedness of
the problem. 
In Section \ref{s3}, 
an operator form of the dual
variational problem without assuming the coercivity condition is considered.
We examine the uniqueness of the solution and describe
all possible multiple solutions for a given right-hand side.
All admissible right-hand sides are determined.
We formulate several equivalent 
necessary and sufficient conditions of wellposedness
in terms of closedness of relevant subspaces. 
We also derive a geometrical condition---a positiveness of a 
\emph{minimum gap} \cite{kato} between relevant closed subspaces. 

A possible application of our theory is the Hellinger--Reissner formulation, 
e.g.,\ \cite{MR1930384}, 
of nonhomogeneous Lam\'e equations for media with (almost) rigid inclusions, 
where the Lagrange multiplier is the displacement, and we get an operator equation
for the stress on the closed subspace of divergence free (in a weak sense) stresses.
Infinitely large Lam\'e coefficients $\lambda$ and $\mu$, in a subdomain, result in
a null-space of the operator in the equation for the stress, 
so the inverse operator does not exist and
the problem is not wellposed in a traditional sense.
Our abstract geometrical condition of generalized wellposedness in this
example is equivalent to a possibility of extension of displacements 
preserving the energy norm of the Lam\'e
operator. It has been proved in \cite{bk94,bkp02} that such an extension is  possible
under some assumptions. We expect that in the
limit case of infinitely large Lam\'e coefficients $\lambda$ and $\mu$ in a
subdomain the pseudoinverse of the operator is bounded, 
which makes the problem wellposed for the stresses in the $L_2$ sense, 
i.e. the $L_2$ norm of the stress is stable even if  
the Lam\'e coefficients are large in a subdomain.
We plan to address this application in the future. 

\section{Abstract symmetric saddle point problems}
\label{s2} 
In this section we essentially follow well known arguments, 
e.g.,\ \cite{bf}, with some simplifications due to the symmetry
of the saddle point problem and our unwillingness to introduce dual spaces. 
Straightforward manipulations, using a pair of complementary 
closed subspaces, allow us, as in \cite{k92,MR1971217}, to formulate separate 
equations for the primary unknown and for the Lagrange multiplier 
of the saddle point problem; see, e.g.,\ survey 
\cite[Sec. 6]{MR2168342} for similar matrix null-space methods. 
We start by formulating and investigating the problem  
using bilinear forms, and then repeat the arguments  
for operator-based formulations that are used 
in the last section of the paper. 

\subsection{Formulations using bilinear forms}
\label{ss2.1} 
Let $\H$ and ${\V}$ be two real Hilbert spaces 
with scalar products and norms denoted by 
$(\cdot ,\cdot )_\H, \|\cdot\|_\H$  
and 
$(\cdot ,\cdot)_\V, \|\cdot\|_\V$
correspondingly. 
Let   
$a(\cdot ,\cdot ): \H \times \H \to \R$ and 
$b(\cdot ,\cdot ): \H \times {\V} \to \R$
be two continuous bilinear forms  
with $a(\cdot ,\cdot )$ symmetric and nonnegative definite. 
We consider the following problem: for a given $g \in \H$
and $f \in \H$ find $\s \in \H$, called 
the ``primary unknown,'' and $u \in {\V}$, called 
the ``Lagrange multiplier,''
such that 
\begin{equation}
\left\{ 
\begin{array}{r}
a(\s , \e) + b(\e , u) = (g, \e)_\H, \ 
\forall \e \in \H, \\ 
b(\s - f, v) = 0, \ \forall v \in {\V}.
\end{array}
\right.  \label{b}
\end{equation}
We place the right-hand side $f$ ``inside'' of the form $b$ as it allows us to take 
$f \in \H,$ not to introduce the dual space ${\V}^{\prime},$ and
makes several statements somewhat simpler.
We call (\ref{b}) a saddle point problem, since 
equations (\ref{b}) are the optimality conditions 
and their solution is a saddle point for the Lagrangian, 
e.g.,\ \cite{bf}, defined by 
$a( \s , \s ) +  2 b( \s -f, u ) - 2 (g, \s )_\H.$

We call a linear manifold, not necessarily closed, a ``subspace" and a closed linear
manifold a ``closed subspace.''
Let us introduce a special notation $\N\subseteq\H$ for the closed subspace, which is 
the null-space of the bilinear form $b(\cdot,\cdot)$ with respect 
to its first argument, i.e. 
$\N=\{\e\in\H: b(\e,v)=0, \ \forall v \in {\V}\}.$
Let us denote  by $\P\equiv\N^\perp\subseteq\H$ the closed subspace which is $\H$-orthogonal (complementary) to $\N$. 
Closed subspaces $\N$ and $\P$ play 
important roles in this paper, so 
let us introduce an $\H$-\emph{orthogonal projector} $P$ on $\H$ such that
$\Ker  (P) = \N$ and $\Im  (P) = \P$
and the complementary projector $P^\perp = I-P$ with 
$\Im  (P^\perp)=\N$  
and
$\Ker  (P^\perp)=\P$, 
where by $\Im(P)$ we denote the range of operator $P$ and,  
with a slight abuse of the notation, by  
$\Ker(P)$ we denote the null-space of operator $P$. 
We assume throughout the paper, unless stated otherwise, that a bounded
operator is defined everywhere on a corresponding space. 
As an orthogonal projector, operator $P:\H\to\H$ is bounded 
$\H$-selfadjoint, $P=P^*$, and satisfies $P=P^2.$   

In the first equation of system (\ref{b}), let us split 
it into two equations, by plugging separately  $\e=P\e\in\P$ and  
$\e=P^\perp\e\in\N$ and using the fact that 
$b(P^\perp\e,u)=0,\  \forall \e \in \H$.  
The second equation in system (\ref{b}) 
has a simple equivalent geometric interpretation:
$\s - f \in \N,$ or $(\s - f,\e)_\H=0, \ \forall \e\in\P$.  
We then rewrite system (\ref{b}) in the following equivalent form: 
\begin{equation}
\left\{ 
\begin{array}{r}
a(\s,\e) + b(\e , u) = (g, \e)_\H, \ \forall \e \in \P, \\ 
a(\s,\e)=(g,\e)_\H, \ \forall \e \in \N, \\ 
(\s - f,\e)_\H=0, \ \forall \e\in\P.
\end{array}
\right.  \label{bb}
\end{equation}

Now we make an important observation that we can treat 
the first line in system (\ref{bb}) as an equation 
for the Lagrange multiplier $u$, given the primary 
unknown $\s$, i.e. 
\begin{equation}
b(\e,u) = (g,\e)_\H - a(\s,\e),  \ \forall \e \in \P. 
\label{beu}
\end{equation}
The last two lines in system (\ref{bb}) 
involve neither the Lagrange multiplier $u$, nor 
the  bilinear form $b$,  
and can be used to determine the primary unknown $\s$:
\begin{equation}
\left\{ 
\begin{array}{r}
a(\s,\e)=(g,\e)_\H, \ \forall \e \in \N, \\ 
(\s - f,\e)_\H=0, \ \forall \e\in\P.
\end{array}
\right.  \label{bes}
\end{equation}
System (\ref{bes}) describes, e.g.,\ \cite{bf}, the optimality conditions 
of the constrained minimization problem 
$\inf\ \{a(\s,\s)-2(g,\s)_\H\},  
\s \in \H:(\s-f,\e)_\H=0, \forall \e\in\P.$

\subsection{Operator-based formulations}
\label{ss2.2} 
In addition to the formulations above involving bilinear forms,  
it is convenient to consider equivalent operator-based formulations.    
We associate with the forms $a$ and $b$ two linear continuous operators 
$A:\H \to \H$ and $B:\H \to {\V}$ defined by 
$(A \s , \e)_{\H} = a(\s , \e), \ 
(B \s , v)_{{\V}} = b(\s , v), \ \forall \e, \s \in \H, \ v \in {\V}.$
In this definition of $A$ and $B$ we follow a slightly simplified, e.g.,\
\cite{MR1971217,MR2002909}, rather than standard \cite{bf}, approach, namely, we do not need dual spaces
 $\H^{\prime}$ and ${\V}^{\prime}.$
Now, we reformulate the main statements 
of subsection \ref{ss2.1} using the 
just defined operators $A$ and $B$.  
The following operator formulation 
\begin{equation}
\left\{ 
\begin{array}{r}
A \s + B^* u = g \ {\mbox{ in }} \ \H, \\ 
B (\s - f) = 0 \ {\mbox{ in }} \ {\V}
\end{array}
\right.  \label{o}
\end{equation}
is equivalent to the original problem (\ref{b}) with the bilinear forms, 
where the adjoint operator 
$B^*:{\V} \to \H$ is defined, as usual, by 
$(\s,B^* v)_{\H} = (B \s , v)_{{\V}}, \ \forall \s \in 
\H, \ v \in {\V}. $
The operator $A$ is selfadjoint and nonnegative definite, 
$A=A^* \geq 0 \ {\mbox{ on }} \ \H $
since it is defined by the symmetric and nonnegative definite form $a.$

We notice that the second equation in system (\ref{o}) 
has the same geometric interpretation 
as in the case of bilinear forms-based system (\ref{b}):
$\s - f \in \Ker (B).$
The   null-space $\Ker (B) \subseteq \H$ and its $\H$-orthogonal 
complement $\overline{\Im  (B^*)} \subseteq \H$ 
have already been denoted by $\N$ and $\P$, correspondingly, 
and introduced together with the 
$\H$-\emph{orthogonal projector} $P$ on $\H$ such that
$\N=\Ker  (P) = \Ker  (B)$ and $\P=\Im  (P) = \overline{\Im  (B^*)}$
in the previous subsection. 

We split the first equation in system (\ref{o}) 
in two orthogonal parts corresponding to $\N$ and $\P$, 
using that $PB^*u =B^*u$ and $P^\perp B^*u=0$, 
since $\Im (B^*) \subseteq \P$.
We replace $B$ with $P$, since they share the same null-space,
in the second equation in system (\ref{o}) to get the following equivalent 
form of system (\ref{bb}):    
\begin{equation}
\left\{ 
\begin{array}{r}
P A \s + B^* u = P g \mbox{ in }  \H, \\ 
P^\perp A \s = P^\perp g \mbox{ in }  \H, \\ 
P (\s - f) = 0 \ {\mbox{ in }} \ {\H}.
\end{array}
\right.  \label{oo}
\end{equation}
We notice that the first line in system (\ref{oo}) is an equation 
for the Lagrange multiplier $u$, given the primary 
unknown $\s$, as in (\ref{beu}), i.e. 
$B^* u = P (g-A\s).$

We next discuss the necessary and sufficient conditions from \cite{bf} of wellposedness
of the problem and make it clear why one can
find weaker necessary and sufficient conditions. 
To simplify our arguments, we take advantage in the rest of the paper 
of the split of the original system into separate equations for the   
Lagrange multiplier $u$ and the primary unknown $\s$ 
that we have described in this section. It is important to 
realize, however, that we have not made 
any substitutions, neither in the solutions $u$ and $\s$, nor in 
the right-hand sides $f$ and $g$. So whatever statements 
we next prove concerning the dependence of the solutions $u$ and $\s$ 
on the right-hand sides $f$ and $g$, these statements 
are equally applicable to both the separate equations
and to the original system in either bilinear form- or 
operator-based context. 

\section{Inf-sup or LBB condition}
\label{sis} 
In this section, we discuss a traditional assumption, being recently
referred to as \emph{Ladygenskaya--Babu\v{s}ka--Brezzi (LBB) condition},
see \citet{MR793735,MR0421106,bf}, that
the range of operator 
$B:\H\to\V,$ denoted by $\Im  (B),$ is closed. 
The closedness of a range of a closed operator is ultimately connected to 
the boundedness of the operator (pseudo-)inverse, e.g.,\ \cite{kato}. 

In our specific situation, operator $B$ is bounded with the closed domain $\H$ 
and, thus, is closed, so its (pseudo-)inverse 
$B^{-1}:{\Im(B) \to  \H / \Ker  (B)}$ is also closed. 
It is necessary to use a factor-space here to define the inverse, 
since the standard operator inverse $B^{-1}:\Im(B)\to\H$ 
does not exist if $\Ker(B)$ is nontrivial.  
We note that $\Ker  (B)$ is closed and that 
the factor-space $\H / \Ker  (B)$ is a Hilbert space, 
as is $\H$. 
In a Hilbert space, a convenient set of representants for the classes 
in the factor-space is simply the corresponding orthogonal 
complement, e.g.,  ${\H} / \Ker (B)$ is isometrically isomorphic to 
$\P=(\Ker (B))^\perp\subseteq\H$, so 
we set $\|\s\|_{{\H} / \Ker (B)}=\|P\s\|_\H.$ 
The subspace $\Im  (B)$  is the domain of 
the closed operator $B^{-1}:{\Im(B) \to  \H / \Ker  (B)}$
therefore, $\Im  (B)$ is closed if and only if 
$B^{-1}:{\Im(B) \to  \H / \Ker  (B)}$ 
is bounded. Closedness of $\Im (B)$ is equivalent to
closedness of $\Im (B^*)$, so all the arguments above 
can be equivalently reformulated for the adjoint operator $B^*$ 
and its (pseudo-)inverse. 

When written in terms of 
inequalities involving the bilinear form $b:$ 
\begin{eqnarray*}
\inf_{ \s \in \H } \sup_{ v \in {\V} } \frac { b(\s, v) } {
\Vert \s \Vert_{\H / \Ker  (B) } \Vert v \Vert_{{\V} } }
&=&\inf_{ \s \in \H } \frac { \Vert B\s\Vert_\V  } {
\Vert \s \Vert_{\H / \Ker  (B) } } \\
&=&\frac{1}{\Vert B^{-1}\Vert_{\Im(B) \to  \H / \Ker  (B)}} 
> 0, 
\end{eqnarray*}
or, equivalently, 
\begin{eqnarray*}
\inf_{ v \in {\V} } \sup_{ \s \in \H } \frac { b(\s, v) } 
{\Vert \s \Vert_{\H } \Vert v \Vert_{{\V} / \Ker  (B^*) }}
&=&
\inf_{ v \in {\V} } \frac { \Vert B^* v\Vert_\H } 
{\Vert v \Vert_{{\V} / \Ker  (B^*) }}\\
&=& 
\frac{1}{\Vert B^{-*}\Vert_{\Im(B^*) \to  \V / \Ker(B^*)}}
> 0, 
\end{eqnarray*}
the LBB condition is also known as the \emph{inf-sup} condition, 
see \citet{MR0421106,bf}, 
where \ ${\V} / \Ker (B^*)$ means the factor-space of ${\V}$ 
with respect to the closed subspace $\Ker  (B^*).$ 
We implicitly assume that 
the arguments in the inf-sup formulas above 
and throughout the paper 
do not make both the numerator and the denominator vanish.
In \citet{MR793735}, the inf-sup condition does not appear to 
be explicitly formulated, instead, closedness of a range of 
the gradient operator is investigated in connection with 
wellposedness of the diffusion equation. 

We note that the induced norms of an operator and its adjoint 
are equal, so both inf-sup expressions above are equal to the 
same constant that we call $c_b$. 
If at least one of the spaces $\H$ or $\V$ is finite dimensional then 
the value  $c_b$ is positive automatically, so 
it becomes important how $c_b$ depends on some parameters, e.g., 
on the dimension.  

Let us mention that in many practical applications 
the space ${\V}$ can be naturally defined such that 
${\Ker} (B^*)= \{ 0 \},$ so 
the latter inf-sup expression of the LBB condition takes the form  
\[
\inf_{ v \in {\V} } \sup_{ \s \in \H } \frac { b(\s, v) } {
\Vert \s \Vert_{\H } \Vert v \Vert_{{\V}}} = c_b > 0, 
\]
which can be most often seen in publications on the subject.
We now contribute our own equivalent formulations of the LBB condition. 
\begin{mylemma}
\label{myLBB} Subspaces $\Im  (B )\subseteq\V$  and 
$\Im  (B B^*)\subseteq\V$ are closed simultaneously.   
Moreover, if either of them is closed we have 
$\Im  (B B^*) = \Im  (B).$
\end{mylemma}
\proof
If $B B^* v=0$ then $(B^* v,B^* v)_{ \H }=0,$ i.e.\ $B^* v=0,$
which proves that
$
\Ker  (B B^*) = \Ker  (B^*). 
$
Taking an orthogonal complement to both parts gives
$
\overline {\Im  (B B^*)} = \overline {\Im  (B)} 
$
as the operator $B B^*$ is selfadjoint. Trivially,
$
\Im  (B B^*) \subseteq \Im  (B).
$
If the range $\Im  (B B^*)$  is closed then
$
\overline {\Im  (B)} = \overline {\Im  (B B^*)} =
\Im  (B B^*)  \subseteq \Im  (B),
$
but clearly ${\Im  (B)}\subseteq \overline {\Im  (B)}$,
which proves closedness of $\Im  (B)=  \Im(BB^*).$

To prove the inverse statement, assuming that ${\Im  (B)}$ is closed, 
we invoke the orthogonal decomposition argument%
\footnote{This proof is suggested by an anonymous referee}:
$\H = \Im(B^*) \oplus (\Im(B^*))^\perp = \Im(B^*) \oplus \Ker(B)$
since  $\Im(B)$ and thus $\Im(B^*)$ are closed. 
Multiplying this equality by $B$ gives 
$\Im(B)= B\H = B(\Im(B^*) \oplus \Ker(B)) = B\Im(B^*) =  \Im(BB^*).$
\endproof
 
We use the previous lemma to introduce 
$(B B^*)^{-1} :\Im(BB^*) \to {\V}/ \Ker(B^*)$ 
in the next Lemma \ref{myLBBc}. 
It is necessary to use the factor-space ${\V}/ \Ker(B^*)$ here, 
since the standard inverse $(B B^*)^{-1} :\Im(BB^*) \to {\V}$ 
does not exist if $\Ker(B^*)$ is nontrivial.
 
\begin{mylemma} \label{myLBBc}
Closedness of $\Im  (B) \subseteq {\V}$ is
equivalent to boundedness of the operator 
$(BB^*)^{-1}:\Im(BB^*)\to {\V}/ \Ker(B^*).$
\end{mylemma} 
\proof
By Lemma \ref{myLBB}, closedness of $\Im  (B )\subseteq\V$ is equivalent to closedness of 
$\Im  (B B^*)\subseteq\V.$ We use several well-known statements 
on closed operators, e.g.,\ \cite{kato}, applied to the operator $BB^*$,
that we have already reviewed in the second paragraph 
of this section for the operator  $B.$
The operator $BB^*$ is bounded and has the closed domain $\V$, 
so the operator is closed and its (pseudo-)inverse 
$(B B^*)^{-1} :\Im  (B B^*) \to {\V}/ \Ker(B^*)$ 
with the domain $\Im  (B B^*)\subseteq\V$ is also closed. 
The domain $\Im  (B B^*)\subseteq\V$ 
of the closed operator $B^{-1}:{\Im(B B^*) \to  \H / \Ker  (B)}$
is closed if and only if  the operator is bounded. 
\endproof

If $\Im  (B)$ is closed then, using Lemmas \ref{myLBB} and \ref{myLBBc}, 
$\Im  (B )=\Im  (B B^*)$  and we can derive the following useful formula
\begin{equation}  \label{P}
P=B^* (B B^*)^{-1} B: \H \to \H.
\end{equation}
Indeed, we first note that $\Im((B B^*)^{-1})\subseteq{\V}/ \Ker(B^*)$ 
is multiplied by $B^*$ in (\ref{P}), so the product 
is independent of the choice of a representant from the 
equivalence class ${\V}/ \Ker(B^*)$ and, thus, is correctly defined. 
Second, righ-hand side of (\ref{P}) is a linear and bounded operator 
as a product of linear and bounded operators. Moreover, it is an 
\emph{orthogonal projector} on $\H$ since it is selfadjoint and 
idempotent, and has the null-space the same as the 
orthoprojector $P$ has.  

If the LBB condition is not satisfied, i.e. $\Im(B)$ is not closed, then
the domain of definition of the operator $B^* (B B^*)^{-1} B$ is the
subspace $\Im  (B^*) \oplus \Ker  (B),$ which is not closed, and
formula (\ref{P}), where $P$ is the orthogonal projector on $\H$ 
with $\Ker(P)=\Ker(B),$ clearly does not hold. 

Let us note that in the case of finite dimensional spaces $\H$ and ${\V}$
the range $\Im  (B )$ is evidently closed, the operator $(B^*)^+=(B B^*)^+ B$ is
the well-known \emph{Moore--Penrose pseudo inverse} of 
$B^*,$ and $P=B^* (B^*)^+$ is the well known formula for 
the orthogonal projector onto the range of $B^*$.

If $\s$ is an exact solution of system (\ref{o}), then $u$ in 
(\ref{o}) can be found from the equation 
$B^* u = - A \s + g \in\Im(B^*).$
If $\s$ is an approximate solution of system (\ref{o}) 
such that the condition $A \s - g\in \Im(B^*),$ 
which is necessary and sufficient for the existence of $u$, 
does not hold, then $u$ can be computed from the 
projected equation $B^* u=P(-A\s+g)\in\P.$
Both the original and the projected equations for $u$ 
are 
wellposed by the LBB assumption, i.e.\ 
$\Im(B^*)=\P$ and 
\[
\| u \|_{{\V} / \Ker  (B^*) } \leq \frac { \| a \| } { c_b } \| \s 
\|_{\H} +
\frac { 1 } { c_b } \| g \|_{\H}. 
\]

Whether the LBB assumption is necessary for wellposedness of
the equation for $u$ depends on if the set of all possible 
right-hand sides $g-A\s$ gives the whole subspace 
$\Im  (B^*),$ see \cite{bf}.  
For example, in a practically important case $g=0$ we have  
$B^* u = - A \s = - P A \s \in \Im(PA)\subseteq\Im(P).$ 
If the latter inclusion is strict, it opens up an opportunity 
for a weaker, compared to the original LBB, 
assumption of wellposedness of the above equation for $u.$  

In the present paper, however, we are concerned with finding 
$\s$, not $u$.
The LBB condition for the bilinear form $b$ 
appears to be of no importance for our results in the next
section where we analyze wellposedness of system 
(\ref{o}) with respect to the $\s$ unknown only, assuming that the $u$
unknown is of no interest, or can be found for a given $\s$
using some \emph{postprocessing}.

\section{Coercivity conditions}
\label{s3} 
\subsection{The standard coercivity condition}
\label{ss31} 
We finally get to the main topic of the paper: an assumption on $A$ which is a condition of
wellposedness of (\ref{o}) with respect to $\s.$
For the reader's convenience, we briefly repeat the 
necessary notation and the system of equations for $\s$ 
to make this section self-consistent.  
Let $\H$ be a real Hilbert space and $P$ be an orthoprojector in $\H$
with a null-space $\Ker(P)=\N$ and 
a range $\Im(P)\equiv{\P}$---we emphasize that the range 
of any orthoprojector in a Hilbert space is closed.  
Let $A$ be a linear and bounded operator such that
$0\leq A^{*}=A$ on $\H$.  
The last two lines in system (\ref{oo}) 
represent an operator form of system (\ref{bes});
they do not involve the Lagrange multiplier $u$ 
or the operator $B$ 
and determine the primary unknown $\s \in \H$:  
\begin{equation}
\left\{ 
\begin{array}{r}
P^\perp ( A \s -g ) = 0 \ {\mbox{ in }} \ \H, \\ 
P (\s - f) = 0 \ {\mbox{ in }} \ \H,
\end{array}
\right.   \label{2.2}
\end{equation}
where $g \in  \H$  and $f \in \H$ are given and $P^{\perp }\equiv I-P.$
We can also replace system (\ref{2.2}) with 
the following equivalent single equation:
\begin{equation}  \label{os}
P^{\perp }A \mid_{\N} \psi=
P^{\perp } g-P^\perp A P f \in \N, \quad \s=\psi+Pf,
\end{equation}
where in (\ref{os}) 
we take a restriction of the operator $P^{\perp }A$ on its invariant
closed subspace $\N$, and  
we  are looking for a solution $\psi \in \N$. Then the necessary and
sufficient condition of wellposedness of problem (\ref{os}) for an 
\emph{arbitrary} $g \in \H$ is, clearly, that the range of 
$P^{\perp }A \mid_{\N}$ is ${\N}.$  
This leads to the traditional assumption, see \cite{bf},   
$a(\s,\s)\geq c_a > 0, \forall\s\in\N, \Vert\s\Vert_{\H}=1$ 
or, in an operator form, $A \geq c_a I$ on  $\N\subseteq\H,$ 
since $A$ is selfadjoint nonnegative. 
Thus, this assumption is also \emph{necessary and sufficient} \cite{b74,bf} for 
wellposedness of
system (\ref{o}) with respect to $\s$ for \emph{an arbitrary} $g \in \H.$
In the rest of the section, we analyze the scenario, where 
$A$ is selfadjoint nonnegative on $\H$, but may be degenerate on $\N$, 
so we impose necessary restrictions on $g \in \H,$ and determine 
a generalized coercivity condition that covers the case of the degeneracy. 

\subsection{Existence, uniqueness, and wellposedness}
\label{ss32}   
Before we investigate the existence and uniqueness of the solution $\s,$
we prove the following technical, but important, lemma.  
\begin{mylemma}\label{l.te1}
Let $P$ be an orthoprojector in $\H$
with a null-space $\Ker(P)=\N$ and 
a range $\Im  (P)\equiv {\P}=\N^\perp,$ and $A$ be a linear and bounded
operator such that
$0\leq A^{*}=A$ on $\H$.  Then 
\begin{equation}
\Ker  (P^{\perp}A) \cap\N= \Ker(A) \cap
\N,  \label{2.3.6}
\end{equation}
\begin{equation}
\{\Ker (P^{\perp }A)\cap \N\}^{\perp }
= \overline{\Im (A)+\P},  \label{2.3.4}
\end{equation}
\begin{equation}
\N=\left(\Ker  (P^{\perp}A) \cap\N\right) \oplus
\overline{P^\perp\Im (A)}.
\label{Nod}
\end{equation}
\end{mylemma}
\proof
We first verify (\ref{2.3.6}).  
It follows from $\Ker (P^{\perp }A) \supseteq \Ker(A)$
that the right-hand side of (\ref{2.3.6}) is included in
the left-hand side. To prove the reverse inclusion, let 
$\varphi \in\N$ and $P^{\perp }A\varphi =0$, then
$0=(P^{\perp}A\varphi,\varphi)=(A\varphi,\varphi)=
\|A^{1/2}\varphi\|^2$
(recall that $A\geq 0$). Then $A^{1/2}\varphi=0$ and $A\varphi=0$.
Therefore, equality (\ref{2.3.6}) holds. 

Equality (\ref{2.3.4}) follows from (\ref{2.3.6}), 
by substituting $\Ker(A)=\F^{\perp}$ and  $\overline{\Im (A)}=\F$
in the well-known simple identity 
$\F^{\perp}\cap\P^{\perp}=(\F+\P)^{\perp}$ 
and noting that 
$
(\overline{\Im (A)}+\P)^{{\perp}{\perp}}=
\overline{\overline{\Im (A)}+\P}=
\overline{\Im (A)+\P}
$
by properties of the closure. 

Finally, to obtain the second term in the orthogonal decomposition (\ref{Nod}) 
of $\N$ we see that by (\ref{2.3.4})
$\{\Ker (P^{\perp }A)\cap \N\}^{\perp }\cap \N 
= \overline{\Im (A)+\P}\cap \N;$ at the same time   
\begin{eqnarray*}
\overline{\Im (A)+\P}\cap \N&=& 
\overline{P^{\perp}\Im (A)+\P}\cap \N\\&=& 
\left(\overline{P^{\perp}\Im (A)}\oplus\P\right)\cap \N=
\overline{P^\perp\Im (A)},
\end{eqnarray*}
which completes the proof of the lemma. 
\endproof


We start with the solution uniqueness. 
\begin{mylemma}\label{l2.0}
Suppose that for some fixed $g\in \H$ and  $f\in \H$ there exists a
solution $\s $ of (\ref{2.2}). Then it is unique provided that 
$
\Ker (A)\cap\N=\{0\};  
$
 otherwise, all possible solutions yield the hyperplane 
$\s +\{\Ker (A) \cap\N\}$ 
and there exists the unique normal (with minimal norm in $\H$) solution of 
(\ref{2.2}) 
that can be also defined as a common element of the above
hyperplane and the closed subspace $ \overline{\Im (A)+\P}$,
which is the set of all normal solutions for all possible $f$ and $g$.
\end{mylemma}
\proof
All solutions of (\ref{2.2}) with $g=f=0$ constitute the closed
subspace
$
\Ker (P^{\perp }A)\cap\N  
$
(may be $0$-dimensional), which by (\ref{2.3.6})
is the same as $\Ker (A)\cap\N$. Hence, all solutions of (\ref{2.2}) with the given 
$g$ and $f,$ provided that there exists at least one solution $\s,$ constitute
the hyperplane $\s +\Ker (A)\cap\N.$ It is
known that each closed hyperplane in a Hilbert space has a unique element with the
minimal norm, i.e. the element that is orthogonal to the directing closed subspace 
$\Ker (A)\cap\N$ of the hyperplane. The orthogonal 
complement to the directing closed subspace is already given by (\ref{2.3.4}). 
\endproof

In the rest of the subsection we use the following 
equation equivalent to (\ref{os}):  
\begin{equation} \label{seq}
(P^\perp  A + P) \s = P^\perp g + P f.  \label{single}
\end{equation}
The assumptions on the right-hand side of the system (\ref{2.2}) 
which ensure the existence of a solution 
are rather standard and follow from (\ref{single}) easily. 
\begin{mylemma}\label{l2.3.3}
 For any $f\in \H$ there exists a solution of 
(\ref{2.2}) if and only if
$
g \in \ia + \P,
$
i.e.
$
P^\perp g + P f \in P^\perp \ia + \P = \ia + \P. 
$
\end{mylemma}
\proof
The subspace (not necessarily closed) 
$ P^\perp \ia + \P$ is simply the range of the operator $P^\perp A + P$
of equation (\ref{single}).
\endproof

The subspace $\ia + \P$ that appears in Lemmas \ref{l2.0}
and \ref{l2.3.3} plays the central role in the following necessary and 
sufficient conditions of wellposedness. 
\begin{mytheorem} \label{t2.0}
The following statements are equivalent: 
\begin{enumerate}
\item The subspace $\ia + \P$ is closed.
\item  The subspace $A \N + \P$ is closed.
\item  The subspace $P^{\perp }\ia$ is closed.
\item  The subspace $P^\perp A \N$ is closed.
\item Problem (\ref{single}) with $f\in \H$ and $g \in \ia + \P$
is well-posed in the factor-space, 
$
\| \s \|_{\H / \{\ka\cap\N\}} \leq c ( \| g \| + \| f \| ), 
$
or, equivalently, 
$\| \s \| \leq c ( \| g \| + \| f \| )$  
for the normal solution $\s \in \overline{\ia + \P}.$
\end{enumerate}
\end{mytheorem}
\proof
(1)$\Leftrightarrow$(3)
We have
$ $\ia + \P$ = P^{\perp }\ia \oplus \P.$\\
(1)$\Leftrightarrow$(2) 
The subspace 
$ P^\perp \ia \oplus \P =  \ia + \P$ is  the range of the operator $P^\perp A + 
P.$
The range of a bounded operator is closed if and only if the
range of the conjugate operator is closed.\\
(2)$\Leftrightarrow$(4)
Using the same arguments as above,
$ A \N + \P = P^\perp A \N \oplus \P.$\\
(1)$\Leftrightarrow$(5)
The operator  $P^\perp A + P$ is bounded and defined everywhere on a
Hilbert space, thus it is closed. Therefore, the (pseudo)inverse operator
$$
(P^\perp A + P)^{-1}: \Im (P^\perp A + P) \rightarrow \H / \{ \ka \cap \N \}
$$
is closed. It is bounded if and only if its domain of definition $ \Im (P^\perp 
A + P)$
is closed. A normal solution is a convenient representant of a factor-class
in a Hilbert space.
\endproof

\subsection{Generalized coercivity conditions}
\label{ss33} 
Statements (1)--(4) in Theorem \ref{t2.0} may not be so 
easily verifiable in practice, so we want to find 
a somewhat easier assumption that generalizes the 
standard coercivity assumption  $A \geq c_a I$ on $\N\subseteq\H,$ 
which itself does not hold if the operator $A$ vanishes on 
a nontrivial subspace of $\N\subseteq\H.$

Let us return back to equation (\ref{os}). We remind the reader that 
the first equation in (\ref{2.2}) is equivalent to the
orthogonal expansion
$\s=\psi+Pf,$ where $\psi=P^\perp\s\in\N.$ 
This and the second equation in (\ref{2.2}) lead to (\ref{os})
that we present here, introducing a special notation 
$K=P^{\perp }A \mid_{\N}$, in the equivalent form 
\begin{equation}
K\psi=\phi,\psi\in \overline{P^\perp \ia},\phi = P^{\perp }g 
-P^\perp A P f\in \overline{P^\perp \ia}.  \label{2.26}
\end{equation}
under the assumption that $g \in \ia+\P.$

The operator $K$ is bounded, selfadjoint, and nonnegative definite on 
$\N$, where $\N\subseteq\H$ inherits the scalar product and the norm of $\H,$
so there exists a bounded, selfadjoint, and nonnegative definite
square root $\sqrt{K}$ on $\N$. Applying the inf-sup condition to the 
operator  $\sqrt{K}$ on $\N$, by direct analogy with 
Lemmas \ref{myLBB} and \ref{myLBBc} and their proofs, we have 
that  $\Ker\left(\sqrt{K}\right)=\Ker(K)$ and 
\begin{mytheorem}
\label{mycc} 
The following statements are equivalent: 
\begin{enumerate}
\item  The subspace $\Im\left(\sqrt{K}\right)\subseteq\N$ is closed.
\item  The subspace $\Im(K)\subseteq\N$ is closed.
\item  The inf-sup condition for the operator  $\sqrt{K}$ on $\N$
\begin{eqnarray}\nonumber
\inf_{ \e \in \N } \sup_{\s \in \N } \frac { \left(\sqrt{K}\e, \s\right)_\N } 
{\Vert \e \Vert_{\N / \Ker (K) } \Vert \s \Vert_{{\N} } }
&=& \inf_{ \e \in \N} \frac { \left\Vert\sqrt{K}\e\right\Vert_\N  } {
\Vert \e \Vert_{\N / \Ker  (K) } } 
\equiv  \frac 1 {\sqrt{\rho}}\\&>& 0 \label{infsupsqrtK}
\end{eqnarray}
holds. 
\item  The norm of the operator $K^{-1} :\Im(K)\to{\N}/ \Ker(K)$ is  
equal to $\rho<\infty.$
\end{enumerate}
Moreover, under either of the assumptions we have 
$\Im\left(\sqrt{K}\right)=\Im(K).$
\end{mytheorem}  

Noticing that $\Im(K)=P^\perp A \N$, we immediately 
see that statements (4) in Theorem \ref{t2.0}
and (2) in Theorem \ref{mycc} are the same, so 
all statements of Theorems \ref{t2.0} and 
\ref{mycc} are equivalent. 
Our last goals in this subsection are 
to present statement (3) of Theorem \ref{mycc} 
in original terms, so that it resembles the 
coercivity condition, and to bound the norm of  
the solution in terms of the norms of the 
right-hand sides, using statement (4) of Theorem \ref{mycc}. 
\begin{mytheorem} \label{t2.1}
For any $g \in \ia + \P$ the following assumption 
\begin{equation}
A\geq \frac 1\rho I \, \mbox{ on the subspace } \, P^{\perp} \ia
\label{2.21}
\end{equation}
with a (finite) constant $\rho >0$
is necessary and sufficient for the normal solution $\s$ 
with $P^{\perp }\s \in P^{\perp } \ia$
to exist and to be unique 
and continuous in $f\in\H$ and $g\in\ia+\P.$ 
Moreover, assumption (\ref{2.21}) implies  
\begin{equation}
\|\s\|^2 \leq \|f\|^2  + \rho^2 \| g - AP f\|^2.
\label{2.22}
\end{equation}
\end{mytheorem}
\proof
First, we note that inequality (\ref{2.21}) on the subspace $P^{\perp } \ia$
is equivalent to the same inequality on 
its closure $\overline{P^{\perp } \ia}$ because of the continuity of $A$
and the scalar product.  
Second, as $(\e,K\e)=(\e,P^\perp A\e)=(\e,A\e)$ for all 
$\e\in\overline{P^\perp\ia}\subseteq\N$, inequality (\ref{2.21})
is also equivalent to 
\begin{equation}
K\geq \frac{1}{\rho}I\mbox{ on the closed subspace }\overline{P^\perp \ia}.
\label{2.25}
\end{equation}

Now we show that (\ref{2.25}) is equivalent to (\ref{infsupsqrtK}), 
which is condition (3) of Theorem \ref{mycc}.  
For the numerator in (\ref{infsupsqrtK}), we have 
$\left\Vert\sqrt{K}\e\right\Vert_\N^2=(\e,K\e).$
To handle the denominator in (\ref{infsupsqrtK}), 
we remind the reader the orthogonal decomposition
$\N=\left(\Ker  (P^{\perp}A) \cap\N\right) \oplus
\overline{P^\perp\Im (A)}$ stated as (\ref{Nod}) 
and proved in Lemma \ref{l.te1}. Splitting $\e\in\N$ 
according to this orthogonal decomposition, we see that its 
first component---from $\Ker(P^{\perp}A)\cap\N=\Ker(K)$---vanishes 
both in the numerator, since it is in the null-space 
of $K$, and in the denominator of (\ref{infsupsqrtK}), 
by the definition of the factor-norm, which gives  
(\ref{2.25}), where only the second component---from 
$\overline{P^\perp \ia}$---survives. 

We conclude that (\ref{2.21}) is equivalent to 
(\ref{infsupsqrtK}), which is condition (3) of Theorem \ref{mycc}, 
and thus, to all statements of Theorems \ref{t2.0} and 
\ref{mycc}. Finally, if (\ref{2.21}) holds then 
the subspace $\Im (P^{\perp }A)$ is closed, the operator 
$K: \overline{P^\perp \ia} \mapsto \overline{P^\perp \ia}$ is an isomorphism
and problem (\ref{2.26}) is wellposed for 
$f\in\H$ and $g\in\ia+\P$, i.e.   
\begin{equation}
\Vert \psi \Vert \leq \rho \Vert P^{\perp } g - P^{\perp }AP f \Vert
\leq \rho \Vert  g - AP f \Vert  \label{2.28}
\end{equation}
by Theorem \ref{mycc}. 
Estimate (\ref{2.22}) follows from $\s=\psi+Pf$
and (\ref{2.28}) due to the statement of Lemma \ref{l2.0}
that the normal solution 
$\s \in \overline{\Im (A)+\P}$, that is, 
$\psi \in \overline{P^{\perp} \ia} = 
P^{\perp }\overline{\Im (A)+\P} =
\Pp \cap ( \Pp \cap \Ker ( A))^\perp$
is the corresponding part of the orthogonal expansion  $\s=\psi+Pf$ for the normal solution. 
\endproof

\subsection{Minimum gap between subspaces}
\label{ss34} 
The rest of the section concerns 
the case where the range of $A$ is closed, so assumption (\ref{2.21})
can be equivalently reformulated using the 
minimum gap between some relevant subspaces. 
We first find a simple way to check if the 
range of $A$ is closed. 
\begin{mylemma} \label{l2.1}
Condition 
\begin{equation}
A \geq \frac 1{\rho _{{\cal D}}}I \, \mbox{on the subspace} \, \Im (A)
\equiv {\bf D}  \label{2.29}
\end{equation}
with a (finite) constant $\rho _{{\cal D}}>0$ is equivalent to closedness of ${\bf D}.$
\end{mylemma}
\proof The operator $A$ is a linear, bounded, and everywhere defined.
Thus, it is closed and its inverse $A^{-1}:\D\to\H/\Ker(A)$
is also closed. 
Boundedness of the inverse is equivalent, on the one hand, to condition 
(\ref{2.29}) and,
on the other hand, to closedness of ${\bf D}.$ 
\endproof

Now we are ready to present a simplified version of the necessary and sufficient 
condition of wellposedness (\ref{2.21}), 
assuming that the range of $A$ is closed. 
\begin{mytheorem} \label{l2.2}
Let the range $\Im (A) \equiv {\bf D}$ be closed, the
orthoprojector on ${\bf D}$ be denoted by $D,$ 
and the constant $\rho _{{\cal D}}>0$ be defined by  (\ref{2.29}). 
Then inequality (\ref{2.21}) is equivalent to the inequality
\begin{equation}
\kappa \equiv
\inf\limits_{\psi \in P^\perp \D 
}
\frac {\| D \psi \| } { \| \psi \| } > 0.
\label{2.30}
\end{equation}
In particular, (\ref{2.29}) and (\ref{2.30})
lead to (\ref{2.21}) with $\rho= \rho_{{\cal D}}/\k^2.$
\end{mytheorem}
\proof
We have $\Im (P^{\perp }A) = P^{\perp }\Im (A) =
P^{\perp}{\bf D},$ i.e. the subspaces indicated in (\ref{2.21}) and 
(\ref{2.30}) coincide. Now the main assertion of the Lemma is a consequence
of relations 
$
\frac {1}{\rho_{{\cal D}}} (D\psi,D\psi)\leq (A\psi,\psi)=
(AD\psi,D\psi)\leq \|A\|(D\psi,D\psi) 
$
which hold for an arbitrary $\psi \in \H$.
\endproof

The next two lemmas provide
alternative assumptions, 
equivalent to (\ref{2.30}),   
which are necessary and sufficient 
for wellposedness, assuming that the range of $A$ is closed. 
It is important to have a choice of a criterion that may be 
easier to check in a practical application. 
For aesthetic reasons we denote $\N\equiv\Pp$.
\begin{mylemma} \label{l2.3} 
Let  $D$ and $P$ be orthogonal projectors onto  closed subspaces 
$\D$ and $\P,$ and  let $D^\perp = I - D$ and $P^\perp = I - P$
be orthogonal projectors onto the orthogonal complements
$\Dp$ and $\Pp$, respectively. 
The following statements are equivalent:
\begin{enumerate}
\item  The subspace $P^{\perp }{\bf D}$ is closed.
\item  The subspace $\D + \P$ is closed.
\item  The subspace $\Dp + \Pp$ is closed.
\item  The subspace $P \Dp$ is closed.
\end{enumerate}
\end{mylemma}
\proof
(1)$\Leftrightarrow$(2) 
The subspace $P^{\perp }{\D}$ is closed iff the subspace
$ P^{\perp }{\D} \oplus \P = {\D} + \P$ is closed
as the terms are orthogonal in the first expression. \\
(2)$\Leftrightarrow$(3)
By Theorem IV-4.8 of \cite{kato}, a sum of closed
subspaces in a Hilbert space is closed if and only if the sum of their
orthogonal complements is closed.\\
(3)$\Leftrightarrow$(4)
Using the same arguments as above,
$\Pp + \Dp = \Pp \oplus P  \Dp.$
\endproof
\begin{mylemma} \label{l2.3c} 
Using the notation of Lemma \ref{l2.3}, 
the following equalities hold: 
\begin{eqnarray*}
\inf\limits_{\psi \in \P , \  \psi \not\in \D } 
\frac {{\rm dist} \{ \psi;\D \} } {{\rm dist}\{\psi;{\D}\cap \P \} }
=
\inf\limits_{\psi \in \D, \  \psi \not\in \P } 
\frac {{\rm dist} \{\psi;\P \} } {{\rm dist}\{\psi;{\P}\cap \D \} }\\
= 
\inf\limits_{\psi \in \Dp, \  \psi \not\in \Pp } 
\frac {{\rm dist} \{\psi;\Pp \} } {{\rm dist}\{\psi;{\Pp}\cap \Dp \} }\\
=
\inf\limits_{\psi \in \Pp, \  \psi \not\in \Dp } 
\frac {{\rm dist} \{ \psi;\Dp \} } {{\rm dist}\{\psi;{\Dp}\cap \Pp \} }\\
=
\inf\limits_{\psi \in P^\perp \D 
}
\frac {\| D \psi \| } { \| \psi \| }
=
\inf\limits_{\psi \in D^\perp \P 
}
\frac {\| P \psi \| } { \| \psi \| }\\
=
\inf\limits_{\psi \in P \Dp 
}
\frac {\| D^\perp \psi \| } { \| \psi \| }
=
\inf\limits_{\psi \in D \Pp 
}
\frac {\| P^\perp \psi \| } { \| \psi \| }. 
\end{eqnarray*}

Moreover, each statement in the previous Lemma is equivalent to 
the positiveness $\k > 0$ in (\ref{2.30}).
\end{mylemma}
\proof
The first three equalities are derived in Section IV-4 of \cite{kato} 
on the \emph{minimum gap} between subspaces, along with a statement that
positiveness of the minimum gap between two given subspaces is a necessary 
and sufficient condition of the sum of the subspaces, 
in our case, $\D + \P$, to be closed. 
We now prove that
$$
\inf\limits_{\psi \in \Pp, \  \psi \not\in \Dp } 
\frac {{\rm dist} \{ \psi;\Dp \} } {{\rm dist}\{\psi;{\Dp}\cap \Pp \} }
=
\inf\limits_{\psi \in P^\perp \D \setminus \{ 0 \}}
\frac {\| D \psi \| } { \| \psi \| }.
$$
All other equalities can be then trivially derived from the previous ones
just by interchanging $P$ and $D.$

We first notice that in the right-hand side we can apply the inf to the closure 
$\overline{ P^\perp \D }  \setminus \{ 0 \}$ as well, because a norm
is a continuous function,
$$
\inf\limits_{\psi \in P^\perp \D \setminus \{ 0 \}}
\frac {\| D \psi \| } { \| \psi \| }
=
\inf\limits_{\psi \in \overline{P^\perp \D} \setminus \{ 0 \}}
\frac {\| D \psi \| } { \| \psi \| }.
$$
We have,
$
\overline{P^\perp \D} = \Pp \cap ( \Pp \cap \Dp )^\perp
$
as
$
\Ker (D \Pp) = \P \oplus ( \Pp \cap \Dp ).
$
The latter can be checked directly.

We always have
$
{\rm dist} \{ \psi;\Dp \} = \| D \psi \|.
$
If 
$
\psi \in \overline{P^\perp \D} =  \Pp \cap ( \Pp \cap \Dp )^\perp 
\subseteq  ( \Pp \cap \Dp )^\perp, 
$
we also have 
$
{\rm dist}\{\psi;{\Dp}\cap \Pp \} = \| \psi \|.
$
Thus,
$$
\frac {{\rm dist} \{ \psi;\Dp \} } {{\rm dist}\{\psi;{\Dp}\cap \Pp \} }
=
\frac {\| D \psi \| } { \| \psi \| }, \, 
\psi \in  \overline{P^\perp \D} \setminus \{ 0 \}.
$$
Finally, using the orthogonal representation
$
\Pp = ( \Pp \cap \Dp ) \oplus  \overline{P^\perp \D},
$
every $\varphi \in \Pp$ can be written as the orthogonal sum
$
\varphi = (\varphi - \psi ) \oplus \psi,$ where 
$\varphi - \psi \in  \Pp \cap \Dp, \,\psi \in \overline{P^\perp \D}.$
Then
${\rm dist} \{ \psi;\Dp \} = {\rm dist} \{ \varphi;\Dp \}$
and also
${\rm dist}\{\psi;{\Dp}\cap \Pp \} = {\rm dist}\{\varphi;{\Dp}\cap \Pp \};$
so the value of the ratio
\[
\frac {{\rm dist}\{\psi ;\Dp\}} {{\rm dist}\{\psi ;\Dp \cap \Pp \} } 
= 
\frac {{\rm dist}\{\varphi;\Dp\} } {{\rm dist}\{\varphi ; \Dp \cap \Pp \}   } 
\]
does not depend on $\varphi -\psi $ and its two infimum values, taken with
respect to 
$
\psi \in \overline{P^\perp \D} \setminus \{ 0 \}
$ 
and 
$\varphi \in \Pp, \, \varphi  \not\in \Dp,$ 
coincide.
\endproof

Finally, we notice that $g=0$ if we apply a saddle point approach
to diffusion or linear elasticity equations. Indeed, 
in the Hellinger--Reissner formulation of
nonhomogeneous Lam\'e equations,  our $\s$ represents the stress tensor,  
the Lagrange multiplier $u$ is the displacement, and 
if we also introduce the stain  $\e$ by the 
stain-displacement relation $\e = - B^* u$, then  
the first line in system (\ref{o})  becomes $A \s - \e =g$, 
which is the constitutive equation (3-D Hooke's law), where of course $g=0$.   
The second line in (\ref{o}) is the equilibrium equation, where 
all body and traction forces are represented by $f\neq0$. 
The assumption $g=0$ allows us to look
for even weaker conditions of wellposedness that we plan to 
investigate in the future.  
\vspace{-1cm}
\section*{Acknowledgments}
The author thanks Ivo Babu\v{s}ka and Franco Brezzi for discussions.
This work has been stimulated by collaboration with Nikolai S. Bakhvalov. 
The author thanks an anonymous referee, who has made numerous useful suggestions to 
improve the original version of the paper, and CU-Denver students Donald McCuan and Christopher Harder for proofreading the paper. 

\begin{thebibliography}{17}
\providecommand{\natexlab}[1]{#1}
\providecommand{\url}[1]{\texttt{#1}}
\expandafter\ifx\csname urlstyle\endcsname\relax
  \providecommand{\doi}[1]{doi: #1}\else
  \providecommand{\doi}{doi: \begingroup \urlstyle{rm}\Url}\fi

\bibitem[Arnold and Winther(2002)]{MR1930384}
D.~N. Arnold and R. Winther.
\newblock Mixed finite elements for elasticity.
\newblock \emph{Numer. Math.}, 92\penalty0 (3):\penalty0 401--419, 2002.

\bibitem[Babu{\v{s}}ka and Aziz(1972)]{MR0421106}
I. Babu{\v{s}}ka and A.~K. Aziz.
\newblock Survey lectures on the mathematical foundations of the finite element
  method.
\newblock In \emph{The mathematical foundations of the finite element method
  with applications to partial differential equations}, pages 1--359. Academic Press, New York,
  1972.
\newblock With the collaboration of G. Fix and R. B. Kellogg.

\bibitem[Bakhvalov and Knyazev(1990)]{bk90}
N.~S. Bakhvalov and A.~V. Knyazev.
\newblock A new iterative algorithm for solving problems of the fictitious flow
  method for elliptic equations.
\newblock \emph{Soviet Math. Doklady}, 41\penalty0 (3):\penalty0 481--485,
  1990.

\bibitem[Bakhvalov and Knyazev(1994{\natexlab{a}})]{bk94}
N.~S. Bakhvalov and A.~V. Knyazev.
\newblock Fictitious domain methods and computation of homogenized properties
  of composites with a periodic structure of essentially different components.
\newblock In Gury~I. Marchuk, editor, \emph{Numerical Methods and
  Applications}, pages 221--276. CRC Press, Boca Raton, 1994{\natexlab{a}}.

\bibitem[Bakhvalov and Knyazev(1994{\natexlab{b}})]{bk95}
N.~S. Bakhvalov and A.~V. Knyazev.
\newblock Preconditioned iterative methods in a subspace for linear algebraic
  equations with large jumps in the coefficients.
\newblock In D.~Keyes and J.~Xu, editors, \emph{Domain Decomposition Methods in
  Science and Engineering}, volume 180 of \emph{Contemporary Mathematics},
  pages 157--162. American Mathematical Society, Providence,
  1994{\natexlab{b}}.

\bibitem[Bakhvalov et~al.(1991)Bakhvalov, Knyazev, and Kobel'kov]{bkk91}
N.~S. Bakhvalov, A.~V. Knyazev, and G.~M. Kobel'kov.
\newblock Iterative methods for solving equations with highly varying
  coefficients.
\newblock In Roland Glowinski, Yuri~A. Kuznetsov, G{\'e}rard~A. Meurant,
  Jacques P{\'e}riaux, and Olof Widlund, editors, \emph{Fourth International
  Symposium on Domain Decomposition Methods for Partial Differential
  Equations}, pages 197--205, Philadelphia, PA, 1991. SIAM.

\bibitem[Bakhvalov et~al.(2002)Bakhvalov, Knyazev, and Parashkevov]{bkp02}
N.~S. Bakhvalov, A.~V. Knyazev, and R.~R. Parashkevov.
\newblock Extension theorems for {Stokes and Lame} equations for nearly
  incompressible media and their applications to numerical solution of problems
  with highly discontinuous coefficients.
\newblock \emph{Numerical Linear Algebra with Applications}, 9\penalty0
  (2):\penalty0 115--139, 2002.

\bibitem[Benzi et~al.(2005)Benzi, Golub, and Liesen]{MR2168342}
M. Benzi, G.~H. Golub, and J. Liesen.
\newblock Numerical solution of saddle point problems.
\newblock \emph{Acta Numer.}, 14:\penalty0 1--137, 2005.

\bibitem[Brezzi(1974)]{b74}
F.~Brezzi.
\newblock On the existence, uniqueness and approximation of saddle point
  problems arising from {L}agrangian multipliers.
\newblock \emph{RAIRO Anal. Numer.}, 2:\penalty0 129--151, 1974.

\bibitem[Brezzi and Fortin(1991)]{bf}
F.~Brezzi and M.~Fortin.
\newblock \emph{Mixed and Hybrid Finite Element Methods}.
\newblock Springer--Verlag, New York, 1991.

\bibitem[Ciarlet et~al.(2003)Ciarlet, Huang, and Zou]{MR2002909}
P.~Ciarlet, Jr., J. Huang, and J. Zou.
\newblock Some observations on generalized saddle-point problems.
\newblock \emph{SIAM J. Matrix Anal. Appl.}, 25\penalty0 (1):\penalty0 224--236, 2003.

\bibitem[Kato(1976)]{kato}
T.~Kato.
\newblock \emph{Perturbation Theory for Linear Operators}.
\newblock Springer--Verlag, New--York, 1976.

\bibitem[Knyazev(2003)]{k03}
A.~V. Knyazev.
\newblock Analysis of transmission problems on {L}ipschitz boundaries in
  stronger norms.
\newblock \emph{J. Numer. Math.}, 11\penalty0 (3):\penalty0 225--234, 2003.

\bibitem[Knyazev(1992)]{k92}
A.~V. Knyazev.
\newblock Iterative solution of {PDE} with strongly varying coefficients:
  algebraic version.
\newblock In R.~Beauwens and P.~de~Groen, editors, \emph{Iterative Methods in
  Linear Algebra}, pages 85--89, Amsterdam, 1992. Elsevier.

\bibitem[Knyazev and Widlund(2003)]{kw03}
A.~V. Knyazev and O.~Widlund.
\newblock {L}avrentiev regularization + {R}itz approximation = uniform finite
  element error estimates for differential equations with rough coefficients.
\newblock \emph{Mathematics of Computation}, 72:\penalty0 17--40, 2003.

\bibitem[Ladyzhenskaya(1985)]{MR793735}
O.~A. Ladyzhenskaya.
\newblock \emph{The boundary value problems of mathematical physics}, volume~49
  of \emph{Applied Mathematical Sciences}.
\newblock Springer-Verlag, New York, 1985.

\bibitem[Xu and Zikatanov(2003)]{MR1971217}
J. Xu and L. Zikatanov.
\newblock Some observations on {B}abu\v ska and {B}rezzi theories.
\newblock \emph{Numer. Math.}, 94\penalty0 (1):\penalty0 195--202, 2003.

\end{thebibliography}
\vspace{-1cm}
\def\cprime{$'$}

\end{document}